\documentstyle[12pt]{article}

\newtheorem{theorem}{Theorem}
\newtheorem{lemma}{Lemma}
\newtheorem{prop}{Proposition}

\def\bsq{\hbox{\vrule width 8pt height 8pt depth 0pt}}

\begin{document}

\title{On quasi-Frobenius semigroup algebras}

\author{B.\,V.\,Novikov (Kharkov, Ukraine)}
\date{}

\maketitle
\begin{abstract}
We define quasi-Frobenius semigroups and find necessary and
sufficient conditions under which a semigroup algebra of a
0-cancellative semigroup is quasi-Frobenius.
\end{abstract}

I.\,S.\,Ponizovskii was first who considered quasi-Frobenius semigroup
algebras over a field for finite semigroups. In \cite{pon} he described
the commutative case completely.  Here we study another class of semigroups
($0$-cancellative). As an example we consider so called modifications of
finite groups \cite{k-n, nov1, nov2}.

\medskip

A semigroup $S$ is called {\it $0$-cancellative} if for all $a,b,c\in S$
from $ac=bc\ne 0$ or $ca=cb\ne 0$ it follows $a=b$. {\bf In what follows
$S$ denotes a finite $0$-cancellative semigroup.} As the next assertion shows,
these semigroups can be called elementary analogously to \cite{pon, okn}.

Let $e$ be an identity of $S$, $H$ the subgroup of invertible elements of $S$
(if $S$ does not contain an identity, we set $H=\emptyset$), $N=S\setminus H$.

\begin{lemma}\label{lem-1}
$N$ is a nilpotent ideal of $S$.
\end{lemma}

{\bf Proof.} If $ab=e$ then $bab=b$, whence ($0$-cancellativity)
$ba=e$. Hence every right inverse element is also left inverse.
This implies that $N$ is an ideal.

Let $a\in S$. Since $S$ is finite, $a^m=a^n$ where, e.\,g., $m<n$.
It follows from $0$-cancellativity that either $a\in H$ or
$a^m=0$. Therefore $N$ is a nil-semigroup, and since $|N|<\infty$,
$N$ is nilpotent \cite{jac}. \bsq

If $A$ is a semigroup or a ring we shall denote by $l_A(B)$ and
$r_A(B)$ the left and right annihilators of a subset $B\subset A$,
respectively.

Analogously with Ring Theory, we call a semigroup $S$ {\it quasi-Frobenius}
if $r_Sl_S(R)=R$ for every right ideal $R$ and $l_Sr_S(L)=L$ for every left
ideal $L$ of $S$.

\begin{lemma}\label{lem-2}
If $S$ is quasi-Frobenius then $S$ contains an identity
(so $H\ne\emptyset$).
\end{lemma}

{\bf Proof.} Suppose that $S$ does not have an identity. Then $S$
is nilpotent by Lemma \ref{lem-1}. Let $S^n=0$,  $S^{n-1}\ne 0$.
Then $0=l_Sr_S(0)=l_S(S)\supset S^{n-1}$ which is impossible. \bsq

Let $F$ be a field, $F_0S$ the contracted semigroup algebra of $S$
(i.\,e. a factor of $FS$ obtained by gluing the zeroes of $S$ and
$FS$). For $A\subset S$ we denote by $F_0A$ the image of the
vector space $FA$ in $F_0S$.

\begin{lemma}\label{lem-3}
If $F_0S$ is quasi-Frobenius then $S$ is quasi-Frobenius too.
\end{lemma}

{\bf Proof.} First show that
\begin{eqnarray}
r_{F_0S}(F_0A)&=&F_0r_S(A)\label{eq-1}\\
l_{F_0S}(F_0A)&=&F_0l_S(A)\label{eq-2}
\end{eqnarray}
for any subset $A\subset S$. Indeed,
$$
r_{F_0S}(F_0A)=r_{F_0S}(A)\supset F_0r_S(A).
$$
Conversely, if $\displaystyle {\sum_{s\in S}}\alpha_ss\in r_{F_0S}(A)$
($\alpha_s\in F$) then $\displaystyle {\sum_{s\in S}}\alpha_sas=0$ for
$a\in A$.  Different non-zero summands in left side of this equality
cannot equal one to another because of $0$-cancellativity. Hence
$$
\sum_{s\in S}\alpha_ss=\sum_{s\in r_S(A)}\alpha_ss \in F_0r_S(A).
$$
The equation (\ref{eq-2}) is proved analogously.

Now for right ideal $R\subset S$ we get from (\ref{eq-1}) and
(\ref{eq-2}):
$$
F_0r_Sl_S(R)=r_{F_0S}[F_0l_S(R)]=r_{F_0S}l_{F_0S}(F_0R)=F_0R.
$$

Further, $F_0A=F_0B$ implies $A\cup 0=B\cup 0$ (here 0 is the zero of $S$).
Since $0\in R\subset r_Sl_S(R)$ we have $r_Sl_S(R)=R$.

Similarly $l_Sr_S(L)=L$. \bsq

Let $S$ contains an identity (i.\,e. $H\ne\emptyset$), $N^n=0\ne N^{n-1}$.
Denote $M(S)=N^{n-1}$. Evidently $S\setminus 0$ is a disjoint union of cosets
of $H$. Since $N$ is an ideal, $HM(S)\subset M(S)$ and $M(S)$ consists of
cosets of $H$ as well.

\begin{lemma}\label{lem-4}
If $S$ is quasi-Frobenius then $M(S)=Ha\cup 0$ for any $a\in M(S)\setminus 0$.
\end{lemma}

{\bf Proof.} Let $a\in M(S)$. Since $Ha\cup 0$ is a left ideal in $S$,
$$
Ha\cup 0= l_Sr_S(Ha\cup 0)= l_S(N)\supset M(S).
$$
So $Ha\cup 0= M(S)$. \bsq

\begin{lemma}\label{lem-5}
If $S$ is quasi-Frobenius then for all $a\in M(S)\setminus 0$ and
$b\in S\setminus 0$ there is an unique element $x$ such that $xb=a$.
\end{lemma}

{\bf Proof.} Uniqueness of $x$ follows from 0-cancellativity.

The assertion is evident for $b\in H$.

Let $b\in N^k\setminus N^{k+1}$. If $k=n-1$ the statement follows
from Lemma \ref{lem-4}. We use decreasing induction on $k$,
supposing that our statement holds for bigger $k$'s .

If $Nb=0$ then $b\in r_S(N)= r_Sl_S(M(S))=M(S)$, i.\,e. $k=n-1$,
the case which has already considered. If $Nb\ne 0$ then $cb\ne 0$
for some $c\in N$. In this case $cb\in N^{k+1}$, so accordingly to
the assumption of induction $xcb=a$ for some $x$. \bsq

Now we are able to prove the main result.

\begin{theorem} \label{th-6}
For a finite $0$-cancellative semigroup $S$ the following
conditions are equivalent:

(i) $S$ is quasi-Frobenius;

(ii) $H\ne\emptyset$ and $M(S)$ is the least non-zero ideal;

(iii) $F_0S$ is Frobenius;

(iv) $F_0S$ is quasi-Frobenius.
\end{theorem}

{\bf Proof.} $1\Longrightarrow 2$ was already proved (see Lemmas
\ref{lem-2} and \ref{lem-4}).

$2\Longrightarrow 3.$ We use Theorem 61.3 from \cite{c-r}. Fix
some $a\in M(S)\setminus 0$ and define a linear function on
$F_0S$:
$$
f(\sum_{s\in S\setminus 0}\alpha_ss)=\alpha_a.
$$
Every element from Ker$f$ is of the form $A=\displaystyle
{\sum_{s\ne a}}\alpha_ss$. Let $\alpha_s\ne 0$. By Lemma
\ref{lem-5} there is $x\in S$ such that $xs=a$. Since $a\ne 0$,
$xt\ne a$ for any $t\ne s$. Hence $xA\not\in {\rm Ker}f$, i.\,e.
Ker$f$ does not contain left (and similarly right) ideals. By
Theorem 61.3 \cite{c-r}, $F_0S$ is Frobenius.

$3\Longrightarrow 4$ is evident.

$4\Longrightarrow 1$ was already proved (Lemma \ref{lem-3}). \bsq

\medskip

\centerline{* * *}

\medskip

A vast variety of examples of finite $0$-cancellative semigroups is given
by modifications of groups \cite{k-n, nov1, nov2}. Remind that
{\it a modification} $G(\ast)$ of a group $G$ is a semigroup on the set
$G^0=G\cup \{0\}$ with an operation $\ast$ such that $x\ast y$ is equal
either to $xy$ or to 0, while
$$
0\ast x=x\ast 0=0\ast 0=0
$$
and the identity of $G$ is the same for the semigroup $G(\ast)$.

In other words, to obtain a modification, one must erase the contents of
some  inputs in the multiplication table of $G$ and insert there zeros so that
the new operation would be associative.

It is clear that modifications are 0-cancellative. Denote as above
by $H$ the subgroup of all invertible elements in $G(\ast)$. By
Lemma \ref{lem-1} its complement $N=G(\ast)\setminus H$ is a
two-sided nilpotent ideal if $G$ is finite. Let $N^n=0\ne
N^{n-1}=M$. We shall describe quasi-Frobenius modifications for
$n\le 3$.

If $n=1$ the semigroup $G(\ast)$ turns into the group $G$ with an
adjoint zero, so its algebra $F_0G(\ast)\cong FG$ is always
quasi-Frobenius.

Let $n=2$. If $G(\ast)$ is quasi-Frobenius then $N=M=Ha\cup 0$.
Therefore when $G$ has a subgroup $H$ of the index 2 we can build
a quasi-Frobenius modification $G(\ast)$ giving its multiplication
by
$$
x\ast y=
\left \{ \begin {array}{ll}
xy & {\rm if}\ x\in H \ {\rm or}\ y\in H,  \\
0      & {\rm otherwise.}
                  \end {array} \right.
$$

Let $n=3$. Fix $a\in M\setminus 0$; then
$M=Ha\cup 0=aH\cup 0$. Let $b\in N\setminus N^2$. By Lemma \ref{lem-5}
for every $h\in H$ there is unique $x\in G(\ast)$ such that $x\ast b=ha$.
Therefore for all $x\in G(\ast)$, $b\in N\setminus N^2$ from $xb\in Ha$
it follows $x\ast b\ne 0$. So the operation $\ast$ must have the form
\begin{equation}\label{eq-3}
x\ast y=
\left \{ \begin {array}{ll}
xy & {\rm if}\ x\in H \ {\rm or}\ y\in H \ {\rm or}\ xy\in Ha,  \\
0,      & {\rm otherwise.}
                  \end {array} \right.
\end{equation}
Hence we have
\begin{prop} Let $G$ be a finite group, $H$ its proper subgroup having
the non-trivial normalizator $N_G(H)\ne H$, $a\in N_G(H)\setminus
H$. Then a modification $G(\ast)$ is given by (\ref{eq-3}) is
quasi-Frobenius.
\end{prop}

{\bf Proof.} We need only to check associativity of $\ast$. This
is equivalent to proving the statement
$$
(x\ast y)\ast z=0 \Longleftrightarrow  x\ast (y\ast z)=0.
$$

If $a,b\in N$ then $a\ast b\in M(G(\ast))$; so $x\ast y\ast z=0$ when
$x,y,z\in N$. Hence it is sufficient to establish associativity only in
the case when one out of elements $x,y,z$ is contained in $H$.

Let, e.\,g., $x\in H$. We have to prove that
\begin{equation}\label{eq-4}
xy\ast z= x(y\ast z).
\end{equation}
However, $xy$ and $y$ belong or do not belong to $H$
simultaneously; also $xyz$ and $yz$ belong or do not belong to
$aH$ simultaneously. Therefore
\begin{eqnarray}
&& xy\ast z=0 \Longleftrightarrow  xy\not\in H \ \& \ z\not\in H \
\& \ xyz\not\in H \nonumber\\
&& \Longleftrightarrow  y\not\in H \ \& \ z\not\in H \ \& \
yz\not\in H\Longleftrightarrow y\ast z=0 \Longleftrightarrow
x(y\ast z)=0\nonumber\bsq
\end{eqnarray}

\medskip

We conclude this note by a remark. Theorem \ref{th-6} and the results of
\cite{pon} are too similar, nevertheless the considered classes of semigroups
(commutative and $0$-cancellative) are rather different. It would be of
interest to find a larger class of semigroups which includes both
above-mentioned classes and yields to some analogue of Theorem \ref{th-6}.

I am thankful to Prof. C. P. Milies for the opportunity to take
part in Conference

\end{document}